\documentclass[a4paper]{amsart}

\usepackage[dvips]{graphicx}

\usepackage{epsf,latexsym}

\begin{document}

\title{Cheeger type inequalities \\for high dimensional simplicial complexes}
\author{Satoshi KAMEI}
\date{}
\address{Department of Liberal Arts, Tokyo University of Technology, 1404-1 Katakura, Hachioji-shi Tokyo 192-0982, Japan}

\email{kamei@stf.teu.ac.jp}


\maketitle

\begin{abstract}
Cheeger inequality is a classical result emerging from the isoperimetric problem in the field of geometry. 
In the graph theory, a discrete version of Cheeger inequality was also studied deeply and the notion was further extended 
for higher dimensional simplicial complexes in various directions.   
In this paper, we consider an analogue of discrete Cheeger inequality for high dimensional simplicial complexes from 
a combinatorial viewpoint. 

\end{abstract}

\section{Introduction}
The isoperimetric problem has long been studied in geometry, and among the results that have been generated within it, Cheeger inequality(\cite{Che}) is one of the most important outcomes. 
Also, the same type inequality on graph, which relates the expansion properties of graphs and the spectra of their Laplacian, was deeply explored, for example, in \cite{Al},\cite{AM},\cite{Dod},\cite{Ta}. 
In this paper, we consider a higher-dimensional analogue of the Cheeger inequality on graph. There have been several studies on higher dimensional versions. 
One approach is considering {\it coboundary expansion}, originating in \cite{Gro}, \cite{LM}, and researched in, for example, \cite{DK}, \cite{SKM}. 
For a combinatorial approach, prominent examples include \cite{PRT}(see also \cite{L}). In \cite{PRT}, a generalization of Cheeger constant was defined as 
\[h(X)=\min_{V=\sqcup_{i=0}^d Ai}\frac{|V| \cdot |F(A_0,A_1,\cdots,A_d)|}{|A_0|\cdot |A_1|\cdot \cdots \cdot |A_d|} \]
for a finite $d$-dimensional simplicial complex $X$ with the set of vertices $V$. Here, the minimum is taken over all partitions of $V$ into nonempty sets $A_0, \cdots, A_d$,
and $F(A_0,\cdots, A_d)$ denotes the set of $d$-dimensional faces with one vertex in each $A_i$. 
The value of $h(X)$ is estimated by the spectral gap of Laplace operator when 
the complex has a complete skeleton in \cite{PRT}, and this result was extended for the case where complexes do not have complete skeletons in \cite{GS}. 

In this work, we consider a modification of $h(X)$ defined above. The reasons for the consideration of taking over a combinatorial approach such as \cite{PRT} are as follows. The first is that the same as \cite{PRT} the area of a region should be measured by the number of vertices, as vertices are often used to represent data in the field of computer science applications. 
The second is that also the same as \cite{PRT} the size of the partitioning of the vertices should be measured by the number of n-dimensional faces. For example in \cite{DK}, there has been research using discrete Hodge theory for ranking among data, where $2$-faces play a significant role. Therefore, it is important to consider the relationship between vertices and $n$-dimensional faces for future applications. Further in this work, we divide 
vertices into two subsets even if we consider the case where the dimensions of simplicial complexes are greater than or equal to $2$, because of  the ease of use in application.

With these in mind, we define a Cheeger type constant for high dimensional simplicial complexes and estimate the values. Let $X$ be a simplicial complex. If 
all inclusion-maximum faces of $X$ have the same dimension, then $X$ is called {\it pure}. 

\vspace*{\baselineskip}
\noindent
\noindent
{\sc Definition 1.1}. For a finite connected pure $n$-simplicial complex $X$, we define 
\[\displaystyle H(X)=\min_{0<|A|<|V|} \frac{|V|\cdot|F(A,V\setminus A)|}{|A|\cdot|V\setminus A|},\]
 where $V$ is the set of all vertices of $X$, $A$ is a nonempty proper subset of $V$ and $F(A,V\setminus A)$ is the set of $n$-faces which contains vertices of  both $A$ and $V \setminus A$. 
\vspace*{\baselineskip}

In the following, $H(X)$ may occasionally be abbreviated as $H$. For estimating the value of $H(X)$, we construct a graph from $X$ as follows.
We set each $(n-1)$-face of $X$ as a vertex and connect these vertices with edges if the corresponding $(n-1)$-faces are contained within an $n$-face.  
We called the graph the {\it embedded graph} of $X$, and the set of the vertices of the embedded graph is denoted as $W$. $\lambda$ denotes the second large eigenvalue
of the adjacent matrix of the embedded graph. We consider only the case where the number of $n$-faces containing an $(n-1)$-faces is fixed and denoted by $D$.  The minimal number of $(n-1)$-faces containing a vertex is denoted by $\delta_{min}$.

For the constant $H$, we prove the following inequalities. First, we consider the case where the dimension of a simplicial complex is $2$. 

\vspace*{\baselineskip}
\noindent
{\bf Theorem 1.2.} {\it Let $X$ be a finite connected pure $2$-dimensional simplicial complex and $V$, $W$, $\delta_{min}$, $D$, $\lambda$ are the same as defined above. Then 
\[\frac{|V|\cdot \delta_{min}\cdot(2D-\lambda)}{4\cdot |W|} \leq H(X).\]}

\vspace*{\baselineskip}

Next, we consider the case where the dimension of a simplicial complex is larger than $2$. Set $\displaystyle k=\left\lfloor (n+1)/2 \right\rfloor$.

\vspace*{\baselineskip}
\noindent
{\bf Theorem 1.3.}  {\it Let $X$ be a finite connected pure $n$-dimensional simplicial complex whose dimension $n$ is larger than $2$ and $V$, $W$, $\delta_{min}$, $k$, $D$, $\lambda$ are the same as defined above. Then
\[\displaystyle \frac{2 \delta_{min }\cdot (nD-\lambda) }{|W|\cdot n \cdot k \cdot(n+1-k)}\leq H(X).\]}

\vspace*{\baselineskip}
To prove Theorem 1.2 and 1.3, we apply the Cheeger inequality on graph to embedded graphs of simplicial complexes. Therefore in section 2, we remind the inequality on graph and prepare 
some notions to consider the higher dimensional cases. We prove Theorem 1.2 in section 3 and Theorem 1.3 in section 4. In section 5, 
we see some examples for the cases where the dimensions of simplicial complexes are $2$.

\section{preliminaries}

Let us start by recalling the Cheeger inequality on graph. Consider a finite connected $d$-regular graph
$G=(V,E)$, and let $A$ be a nonempty proper subset of $V$. We define a quantity called the Cheegar constant as follows.

\vspace*{\baselineskip}
\noindent
{\sc Definition 2.1}. $\displaystyle h(G)= \min_{0<|A|<|V|} \frac{|V|\cdot| E\left(A,V\setminus A\right)|}{|A|\cdot |V\setminus A|}$.

\vspace*{\baselineskip}
This constant satisfies the following inequality, where $\lambda$ denotes the second greater eigenvalue of the adjacent matrix of $G$.

\vspace*{\baselineskip}
\noindent
{\bf Theorem 2.2(\cite{Al},\cite{AM},\cite{Dod}).} 
$\displaystyle d-\lambda \leq h(G) \leq 2\sqrt{2d(d-\lambda)}$.

\vspace*{\baselineskip}

This theorem was proved by Dodziuk in \cite{Dod}, and independently by Alon-Milman in \cite{AM}, and Alon in \cite{Al}.
The proof of this theorem can be found in, for example, \cite{HLW}. Note that in~\cite{HLW}, a slightly different constant $\phi(G)$ is defined 
instead of $h(G)$ as follows: 
\[\displaystyle \phi(G)= \min_{0<|A|<|V|/2} \dfrac{E\left(A,V\setminus A\right)}{|A|}.\] 
It is then shown that the inequality $\phi(G) \leq h(G) \leq 2\phi(G)$ holds, and that 
$\phi(G)$ satisfies $\displaystyle (d-\lambda)/2 \leq \phi(G) \leq \sqrt{2d(d-\lambda)}$.

We introduce 
some notions related to simplicial complexes.  
An {\it (abstract) simplicial complex} $X$ with vertex set $V$ is a collection of subsets of $V$,
called {\it faces} or {\it simplices}, which is closed under taking subsets. That is, if $\sigma \in X$ and $\tau \subset \sigma$,
then $\tau \in X$. The {\it dimension} of a simplex $\sigma$ is $|\sigma|-1$. The subcomplex of $X$ formed by all of $j$-dimensional simplices and their faces  
is called {\it $j$-skeleton} of $X$. The $0$-dimensional simplices in $X$ are {\it vertices} and the $1$-dimensional simplices are {\it edges}. %The inclusion-maximum faces are called {\it facets}. 
The {\it degree} of a $j$-face is the number of $(j+1)$-faces that contain it. For a $j$-face $\sigma$, we denote the degree of $\sigma$ as $\deg(\sigma)$. 
Further in this paper, the number of $(n-1)$-faces that contain a vertex $v$ is denoted as $\delta(v)$. Note that for the case where the dimension of a simplicial complex is $2$ and $v$ is a vertex of the simplicial complex, $\deg(v)$ and $\delta(v)$ determine the same value.  
As we defined in section~1, $\delta_{min} = \min_{v \in V} \delta(v)$, where $V$ is the set of all vertices of a simplicial complex $X$.  
Also as mentioned in section 1, we only consider the case where the degrees of the $(n-1)$-faces of a simplicial complex are constant $D$. Thus the degree of each vertex of the embedded graph is $nD$. For each edge of the embedded graph which  connects two $(n-1)$-faces of $X$,  there is an $n$-face which contains both of the $(n-1)$-faces. 
In this situation, we will refer to this as the edge of the embedded graph being {\it contained} within the $n$-face.
%In this situation, we will call that the edge of the embedded graph is {\it contained} in the $n$-face.  
The embedded graph captures the relationship between 
$(n-1)$-faces and $n$-faces in a simplicial complex. 
Thus we will use the embedded graph to estimate the constant $H$ of a simplicial complex.

\section{An inequality for $2$-dimensional simplicial complexes}

In this section, we estimate $H$ of $2$-dimensional simplicial complexes. 
Let $ G=(V,E)$ be a graph and $V'$ be a subset of $V$. The subgraph of G {\it induced} by $V'$ is the subgraph $G'=(V',E')$ 
such that $G'$ has the set of vertices $V'$ and for all $u$,$v \in V'$, $e=uv \in E'$ if and only if $e \in E$. We call a {\it induced subgraph} of $V'$ 
if the subgraph is induced by $V'$.

\vspace*{\baselineskip}
\noindent
{\bf Lemma 3.1.} {\it Let $X$ be a finite connected pure $2$-dimensional simplicial complex, and let $V$ be the set of all vertices of $X$. 
Consider the $1$-skeleton of $X$ as a graph, denoted by $G_V$. If $A$ is a subset of $V$ that realizes 
$\displaystyle H=\min_{0 < |A|<|V|} \frac{|V|\cdot |F(A,V\setminus A)|}{|A|\cdot |V\setminus A|}$, then the subgraph of $G_V$ induced by $A$ is connected.}

\vspace*{\baselineskip}

\noindent
{\it Proof.} Let $G_A$ denote the subgraph of $G_V$ induced by $A$ and we assume that $G_A$ is not connected. We can split $A$ into two disjoint subsets $A_1$ and $A_2$, where $A =A_1 \cup A_2$ 
and $G_A$ is the disjoint union of the induced subgraphs by $A_1$ and $A_2$. 
Since there are no edges of $G_V$ connecting a vertex of $A_1$ to a vertex of $A_2$, we have 
$F(A,V\setminus A)=F(A_1, V\setminus A_1)\cup F(A_2,V\setminus A_2)$ and $F(A_1, V\setminus A_1)\cap F(A_2,V\setminus A_2)=\emptyset$. 
Let $f_1=|F(A_1, V\setminus A_1)|$, $f_2=|F(A_2,V\setminus A_2)|$, $v=|V|$, $a_1=|A_1|$, $a_2=|A_2|$. 
We assume that $\displaystyle \frac{f_2}{a_2(v-a_2)}\geq \frac{f_1}{a_1(v-a_1)}$, and if this is not the case, we can simply switch $A_1$ and $A_2$. 
Using the assumption, we get $f_2 v a_1 -f_2 a_1^2\geq f_1v a_2-f_1 a_2^2$. Therefore, we have
\[ \frac{|F(A,V\setminus A)|}{|A|\cdot |V\setminus A|}-\frac{ |F(A_1,V\setminus A_1)|}{|A_1|\cdot |V\setminus A_1|}=\frac{f_1+f_2}{(a_1+a_2)(v-a_1-a_2)}-\frac{ f_1}{a_1(v-a_1)}\]
\[\displaystyle =\frac{f_2v a_1-f_2a_1^2+2f_1 a_1 a_2-f_1va_2+f_1a_2^2}{(a_1+a_2)(v-a_1-a_2)a_1(v-a_1)}\]
\[\geq\frac{f_1v a_2-f_1a_2^2+2f_1 a_1 a_2-f_1va_2+f_1a_2^2}{(a_1+a_2)(v-a_1-a_2)a_1(v-a_1)}\]
\[\displaystyle =\frac{2f_1 a_1 a_2}{(a_1+a_2)(v-a_1-a_2)a_1(v-a_1)}>0.\]
This means that 
\[\displaystyle \frac{|V|\cdot |F(A,V\setminus A)|}{|A|\cdot |V\setminus A|}>\frac{|V|\cdot |F(A_1,V\setminus A_1)|}{|A_1|\cdot |V\setminus A_1|}\]
which contradicts the choice of $A$. $\Box$

\vspace*{\baselineskip}

%Notice that we can adapt the statement for $V\setminus A \subset V$. Thus 

The next lemma establishes a correspondence between $2$-dimensional simplicial complexes and their embedded graphs.

\vspace*{\baselineskip}
\noindent
{\bf Lemma 3.2.} {\it Let $X$ be a finite connected pure $2$-dimensional simplicial complex and let $V$ be the set of all vertices of $X$. 
Assume that a nonempty proper subset $A$ of $V$ realizes 
\[\displaystyle H=\min_{0 < |A|<|V|} \frac{|V|\cdot |F(A,V\setminus A)|}{|A|\cdot |V\setminus A|}.\] 
Then there exists a subset $B$ of $W$ satisfying $0<|B|<|W|$ and the following conditions:
\[ \frac{|E(B,W\setminus B)|}{2}\leq |F(A,V\setminus A)| , \ |A|\leq |B| \ ,  |V\setminus A|\leq \frac{2| W \setminus B|}{\delta_{min}}.\]}

\vspace*{\baselineskip}
\noindent
{\it Proof.} % We assume that $|A|\leq|V\setminus A|$, otherwise we can simply swap $A$ and $V\setminus A$. 
For each edge of $X$, we treat it as a vertex of the embedded graph and we allocate the vertex to $B$ or $W\setminus B$ as follows. 
If both end vertices of an edge of $X$ are contained in $A$, we set the edge as a vertex of $B$. 
Similarly, if both end vertices of an edge of $X$ are contained in $V\setminus A$, we set the edge 
as a vertex of $W\setminus B$.
If neither of these cases applies, we choose one edge at random 
and set the edge as a vertex of $B$ and set the remaining edges as vertices of $W\setminus B$. 
%Based on the assumption for $A$ and Lemma 3.1, the subgraph induced by the vertices $B$ 
Based on Lemma 3.1, the subgraph induced by the vertices $B$ 
in the embedded graph of $X$ is connected. Thus $| A| \leq | B|$ and the equality is achieved when the induced subgraph is a tree. 
To deduce that $| V\setminus A| \leq \frac{2 | W\setminus B|}{\delta_{min}}$,  we must compare $2|W\setminus B|$ and $\sum_{v \in V\setminus A} \delta (v)$.
Consider a pair $(v,e)$ where $v \in V$ and $v$ is incident to an edge $e$ of $X$; we refer to this pair as an {\it incident pair}. It's important to note that the count of incident pairs $(v,e)$ with the condition $v \in V\setminus A$ is equal to $\sum_{v \in V\setminus A} \delta (v)$. 
All the edges of $X$ connecting vertices in $V\setminus A$ are contained in $W\setminus B$.
If all the edges connecting vertices in both $V\setminus A$ and $A$ are also contained in $W\setminus B$, every incident pair $(v,e)$ with $v \in V\setminus A$ contributes to $2| W\setminus B|$. Consequently, we can establish the inequality $2| W\setminus B| > \sum_{v \in V\setminus A} \delta (v)$.  
However, there exists exactly one edge connecting a vertex in $V\setminus A$ 
and a vertex in $ A$, which is contained in $B$. Consequently, one incident pair $(v,e)$ with $v \in V\setminus A$ does not contribute to $2| W\setminus B| $. Let's consider a $2$-face of $X$ that contains the particular edge. 
This $2$-face also contains at least one other edge in $W\setminus B$
that connects a vertex in $V\setminus A$ and a vertex in $ A$. See Figure 1. The incident pair involving the vertex in $A$ and this edge contributes to $2| W\setminus B|$ and compensates for the incident pair involving the vertex of $V\setminus A$ from the previous edge. 
As a result, we have $2|W\setminus B|\ge \sum_{v \in V\setminus A} \delta (v)\ge \delta_{\rm min} \cdot | V\setminus A|$, which implies $|V\setminus A| \leq \frac{2 |W\setminus B|}{\delta_{min}}$. 
To prove that $\frac{|E(B,W\setminus B)|}{2}\leq |F(A,V\setminus A)|$, 
let $F'$ be the set of $2$-faces of $X$ that contains edges of $E(B,W\setminus B)$, then each $2$-face of $F'$ contains both a vetex of $V \setminus A$ and a vertex of $ A$ and two edges of $E(B,W \setminus B)$. Thus $F' \subset F(A,V\setminus A)$. This implies $|F(A, V\setminus A)| \geq |F'| = \frac{|E(B,W\setminus B)|}{2}$.~$\Box$
 
%because there may be a $2$-face such that all barycenters of edges of the $2$-face belong to $B$ and two of all vertices belong to $A$ and one belong to $V \setminus A$.  
%Furthermore, each $2$-face of $F'$ contains $2$ edges contained in $E(B,W\setminus B)$, 

 \begin{figure}
   \centerline{ 
    \epsfysize=2.3cm
   \epsfbox{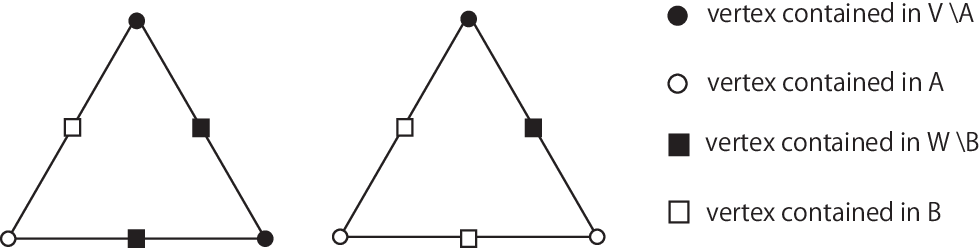}}
   \caption{edges connecting vertices in $A$ and in $V\setminus A$}   
\end{figure}

\vspace*{\baselineskip}

\noindent
{\it Proof of Theorem 1.2.} Let $A\subset V$ realizes $\min_{0<|A|<|V|} \frac{|V|\cdot F(A,V\setminus A)|}{|A|\cdot |V\setminus A|}$. There exists $B\subset W$ which satisfies the conditions of Lemma 3.2. Thus from Lemma 3.2 and Theorem 2.2, we can establish the inequality     
\[H \geq \frac{|V| \cdot \frac{|E(V,W\setminus B)|}{2}}{|B|\cdot \frac{2|W\setminus B|}{\delta_{min}}}=\frac{|V|\cdot \delta_{min}}{4 |W|}\cdot \frac{|W|\cdot |E(B,W\setminus B)|}{|B|\cdot |W\setminus B|}\geq \frac{|V|\cdot \delta_{min}\cdot (2D-\lambda)}{4\cdot |W|}. \ \Box \]

\vspace*{\baselineskip}
\section{An inequality for high dimensional simplicial complexes}

In this section, we discuss the case where the dimensions of simplicial complexes are greater than $2$.

\vspace*{\baselineskip}
\noindent
{\bf Claim 4.1.} {\it If we divide all vertices of an $n$-simplex into two nonempty sets $P$ and $Q$, then the number of $(n-1)$-faces of the $n$-simplex containing vertices of both P and Q is either $n$ or $n+1$.}

\vspace*{\baselineskip}

\noindent
{\it Proof.} If $P$ (or $Q$) has $n$ vertices, then there exists exactly one $(n-1)$-face whose vertices are all contained in $P$ (or $Q$), hence the number of  $(n-1)$-faces containing vertices belonging to both P and Q is $n$.  
In the other cases, all of the $(n-1)$-faces contain vertices belonging to both $P$ and $Q$.$\Box$

\vspace*{\baselineskip}

We prepare a key lemma for the proof of Theorem 1.3. Remind that we define $ k=\left\lfloor (n+1)/2 \right\rfloor$ in section 1. 

\vspace*{\baselineskip}
\noindent
{\bf Lemma 4.2.} Let $X$ be a finite connected pure $n$-simplicial complex and let $V$ be the set of all vertices of $X$. 
Suppose that $A$ is a subset of $V$ satisfying $0 <|A|<|V|$. 
Let $W$ be the set of all vertices of the embedded graph of $X$. Then, there exists a subset $B$ of $W$ that satisfies the 
following conditions:  
\[|E(B,W\setminus B)| \leq k(n+1-k)|F(A,V\setminus A)| ,\ 1 \leq |B|, \ |V\setminus A| \leq \frac{n|W\setminus B|}{\delta_{min}}.\]

\vspace*{\baselineskip}
\noindent
\noindent
{\it Proof.} For each $(n-1)$-face of $X$, we treat it as a vertex of the embedded graph and we allocate the vertex to $B$ or $W\setminus B$ as follows. If all the vertices of an $(n-1)$-face of $X$ are contained in $A$, 
we set the $(n-1)$-face as a vertex of $B$. 
Similarly if all the vertices of an $(n-1)$-face of $X$ are contained in $V\setminus A$, 
we set the $(n-1)$-face as a vertex of $W \setminus B$. 
If neither of these cases applies, we choose one $(n-1)$-face at random 
and set it as a vertex of $B$, and set the remaining $(n-1)$-faces as vertices of $ W\setminus B$. 
Let $F'$ be the set of $n$-faces such that each element of $F'$ contains at least one element of $E(B,W\setminus B)$. Then each $n$-face of $F'$ contains a vertex of $A$ and a vertex of $V\setminus A$. Thus we have $F' \subset F(A, V\setminus A)$. 
%There may be an $n$-face such that some $(n-1)$-faces of the $n$-face contain vertices of both $A$ and $V \setminus A$ and barycenters of all the $(n-1)$-faces of the $n$-face are contained in $B$.  
If an $n$-face of $X$ contains vertices of both $B$ and $W \setminus B$, the number of edges of $E(B,W\setminus B)$ contained in the $n$-face is 
at least $n$ and at most $k(n+1-k)$. Thus, we have $|E(B,W\setminus B)|\leq k(n+1-k)|F'|\leq k(n+1-k)|F(A,V\setminus A)|$. 
To deduce that $|V\setminus A| \leq \frac{n|W\setminus B|}{\delta_{min}}$, we must compare $n|W\setminus B|$ and $\sum_{v \in V\setminus A} \delta(v)$.
Consider a pair $(v,f)$ where $v \in V$ and $v$ is contained in an $(n-1)$-face $f$ of $X$; we refer to this pair as a {\it contained pair}. It's important to note that the count of contained pairs $(v,f)$ with the condition $v \in V\setminus A$ is equal to $\sum_{v \in V\setminus A} \delta (v)$. 
All the $(n-1)$-faces of $X$ containing vertices of $V\setminus A$ are contained in $W\setminus B$. 
If all the $(n-1)$-faces containing vertices of both $V\setminus A$ and  $A$ are also contained in $W\setminus B$, every contained pair $(v,f)$ with $v \in V\setminus A$ contributes to $n| W\setminus B|$. Consequently, $n|W\setminus B|> \sum_{v \in V\setminus A} \delta (v)$. 
However, among the $(n-1)$-faces in $B$, exactly one $(n-1)$-face contains vertices of both $V\setminus A$ and $ A$. Consequently, at most $n-1$ contained pairs containing vertices in $V\setminus A$ do not contribute to $n| W\setminus B| $. 
Consider an $n$-face which contains the particular $(n-1)$-face. From Claim 4.1, 
the number of $(n-1)$-faces of the $n$-face containing vertices in both $A$ and $V\setminus A$ is at least $n$. 
%there are at least $n$ $(n-1)$-faces of the $n$-face containing vertices in both $A$ and $V\setminus A$. 
Therefore, at least $n-1$ of these $(n-1)$-faces are in $W\setminus B$. The contained pairs involving the vertices in $A$ and these $(n-1)$-faces contribute to $n| W\setminus B|$ and compensate for the $n-1$ contained pairs involving the vertices of $V\setminus A$ from the previous $(n-1)$-face. 
As a result, we have $n|W\setminus B| \ge \sum_{v \in V\setminus A} \delta(v) \ge {\delta_{min}} \cdot |V\setminus A|$, which implies $|V\setminus A| \leq \frac{n|W\setminus B|}{\delta_{min}}$. The inequality $| B|\geq 1$ is obvious. 
$\Box$

\vspace*{\baselineskip}

\noindent
{\it Proof of Theorem 1.3.}  We assume that $|A|\leq|V\setminus A|$, otherwise we can simply swap $A$ and $V\setminus A$. Then, 
\[\displaystyle H(X)=\min_{0<|A|<|V|} \frac{|V|\cdot|F(A,V\setminus A)|}{|A|\cdot|V\setminus A|}\]
\[=\min_{0<|A|<|V|} \frac{| F(A,V\setminus A)|}{\frac{|A|}{|V|}\cdot|V\setminus A|}\geq\min_{0<|A|<|V|} \frac{2 |F(A,V\setminus A)|}{|V\setminus A|}.\] 
From Lemme 4.2 and Theorem 2.2, 
\[\min_{0<|A|<|V|} \frac{2 |F(A,V\setminus A)|}{ |V\setminus A|}\geq \frac{2\cdot\frac{|E(B,W\setminus B)|}{k(n+1-k)}}{\frac{n|W\setminus B|}{\delta_{min}}}\]
\[\geq\frac{2\delta_{min }}{|W|\cdot n\cdot k(n+1-k)}\cdot \frac{|W|\cdot |E(B,W\setminus B)|}{|B|\cdot |W\setminus B|}\]
\[\geq \frac{2\delta_{min }\cdot (nD-\lambda) }{|W|\cdot n \cdot k \cdot(n+1-k)}. \Box\]

\vspace*{\baselineskip}

\section{Examples}

In this section, we will be looking at $2$-dimensional simplicial complexes. On the left-hand side of each figure,
we have a simplicial complex, while on the right-hand side, we have its corresponding embedded graph. The black circles 
represent the vertices of the simplicial complex, while the black rectangles 
represent the vertices of the embedded graph.

\vspace*{\baselineskip}

\vspace*{\baselineskip}
 \begin{figure}
   \centerline{ 
    \epsfysize=2.5cm
   \epsfbox{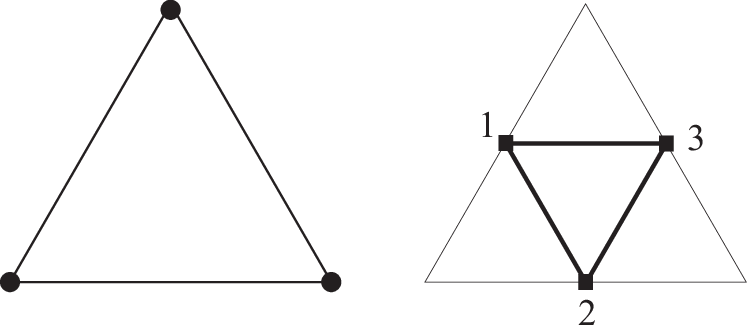}}
   \caption{the simplicial complex in Example 4.1 and its embedded graph}
\end{figure}

\noindent
{\bf Example 4.1.} Let's take a look at a simple yet crucial example. Consider a simplicial complex consisting of only one $2$-simplex. as shown in  
Figure 2. 
It is easy to see that $H= \frac{3\cdot 1}{1\cdot 2}=\frac{3}{2}$. On the other hand, 
the adjacency matrix of the embedded graph is 
$\left( \begin{array}{ccc} 0&1&1 \\ 1 & 0 &1 \\ 1 & 1&0\end{array}\right)$, and the second eigenvalue is $-1$. Therefore  
$\frac{|V|\cdot \delta_{min}}{4 |W|}(2D-\lambda)=\frac{3 \cdot 2}{4 \cdot 3} (2\cdot 1 -(-1))=\frac{3}{2}$.
This example achieves 
the equality of Theorem 1.2.

 \begin{figure}
   \centerline{ 
    \epsfysize=4.5cm
   \epsfbox{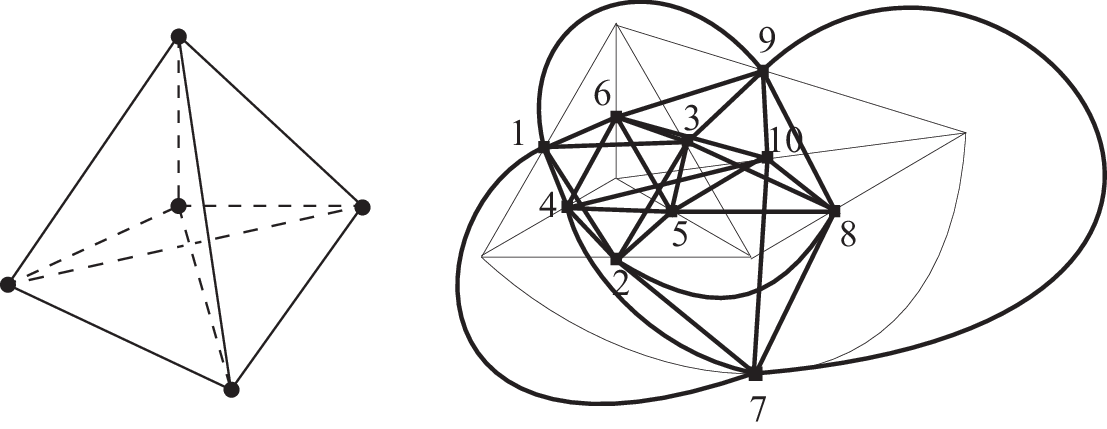}}
   \caption{the simplicial complex in Example 4.2 and its embedded graph}
\end{figure}

\noindent
{\bf Example 4.2.} In the second example, each degree of edge is $3$, as shown in Figure 3. $H= \frac{5\cdot 6}{2\cdot 3}=5$. On the other hand,
the adjacent matrix of the embedded graph is 
\[\left( \begin{array}{cccccccccc} 0&1&1 &1&0&1&1&0&1&0\\
 1 & 0 &1&1&1&0&1&1&0&0 \\
 1 & 1&0&0&1&1&0&1&1&0 \\
1 & 1&0&0&1&1&1&0&0&1\\
0&1&1&1&0&1&0&1&0&1\\
1&0&1&1&1&0&0&0&1&1\\
1&1&0&1&0&0&0&1&1&1\\
0&1&1&0&1&0&1&0&1&1\\
1&0&1&0&0&1&1&1&0&1\\
0&0&0&1&1&1&1&1&1&0
\end{array}\right)\]
and the second eigenvalue is $1$. Therefore $\frac{|V|\cdot \delta_{min}}{4 |W|}(2D-\lambda)=\frac{4 \cdot 4}{4 \cdot 10} (2\cdot 3 -1)=2$.

\vspace*{\baselineskip}
\noindent
{\bf Example 4.3.} The third example is shown in Figure 4. $H= \frac{8\cdot 6}{4\cdot 4}=3$. The adjacent matrix of the embedded graph is
\[\left( \begin{array}{cccccccccccccccccc} 
0&1&1 &1&1&0&0&0&0&0&0&0&0&0&0&0&0&0\\
1&0&1 &0&1&1&0&0&0&0&0&0&0&0&0&0&0&0\\
1&1&0 &1&0&1&0&0&0&0&0&0&0&0&0&0&0&0\\
1&0&1 &0&0&0&1&1&0&0&0&0&0&0&0&0&0&0\\
1&1&0 &0&0&0&0&0&1&1&0&0&0&0&0&0&0&0\\
0&1&1 &0&0&0&0&0&0&0&1&1&0&0&0&0&0&0\\
0&0&0 &1&0&0&0&1&0&0&0&1&1&0&0&0&0&0\\
0&0&0 &1&0&0&1&0&1&0&0&0&0&1&0&0&0&0\\
0&0&0 &0&1&0&0&1&0&1&0&0&0&1&0&0&0&0\\
0&0&0 &0&1&0&0&0&1&0&1&0&0&0&1&0&0&0\\
0&0&0 &0&0&1&0&0&0&1&0&1&0&0&1&0&0&0\\
0&0&0 &0&0&1&1&0&0&0&1&0&1&0&0&0&0&0\\
0&0&0 &0&0&0&1&0&0&0&0&1&0&0&0&1&1&0\\
0&0&0 &0&0&0&0&1&1&0&0&0&0&0&0&0&1&1\\
0&0&0 &0&0&0&0&0&0&1&1&0&0&0&0&1&0&1\\
0&0&0 &0&0&0&0&0&0&0&0&0&1&0&1&0&1&1\\
0&0&0 &0&0&0&0&0&0&0&0&0&1&1&0&1&0&1\\
0&0&0 &0&0&0&0&0&0&0&0&0&0&1&1&1&1&0\\
\end{array}\right)\]

and the second eigenvalue is $1+\sqrt{5}$. Therefore $\frac{|V|\cdot \delta_{min}}{4 |W|}(2D-\lambda)=\frac{8 \cdot 3}{4 \cdot 18} (2\cdot 2 -(1+\sqrt{5}))=\frac{3-\sqrt{5}}{3}$.

 \begin{figure}
   \centerline{ 
    \epsfysize=5cm
   \epsfbox{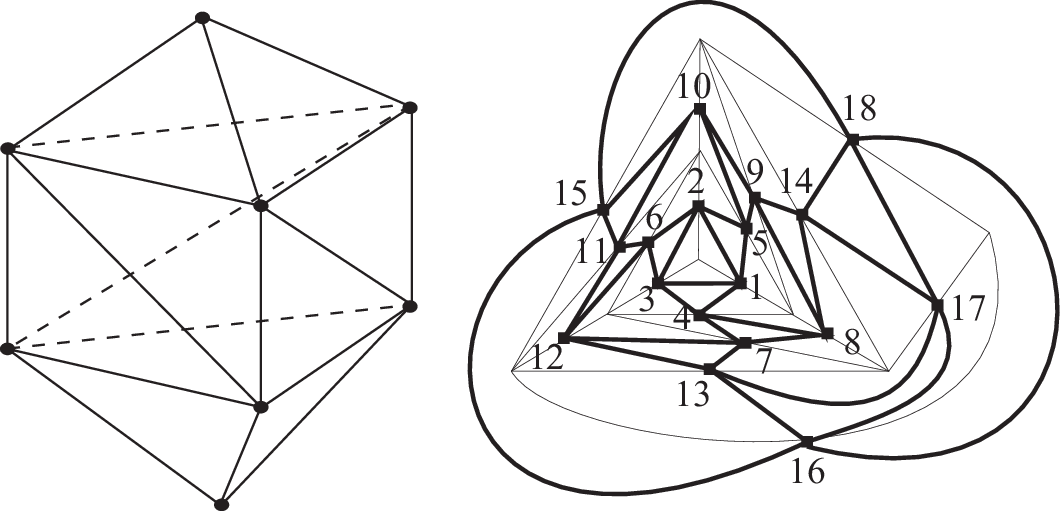}}
   \caption{the simplicial complex in Example 4.3 and its embedded graph}
\end{figure}


\begin{thebibliography}{20}
  \bibitem{Al}
	 N. Alon, {\it Eigenvalues and expanders}, Combinatorica {\bf 6} (1986), no. 2, 83-96.
  \bibitem{AM}
 		N. Alon and V.D. Milman, {\it $\lambda_1$, isoperimetric inequalities for graphs, and superconcentrators}, Journal of Combinatorial Theory, Series B {\bf 38} (1985), no. 1, 73-88.
\bibitem{Che}
  J. Cheeger, {\it A lower bound for the smallest eigenvalue of the Laplacian
Problems in Analysis} Problems in analysis, {\bf 195}(1970), 199.
  \bibitem{Dod}
    J. Dodziuk, {\it Difference equations, isoperimetric inequality and transience of certain random walks}, Trans. Amer. Math. Soc {\bf 284} (1984).
\bibitem{DK}
 D. Dotterrer and M. Kahle, {\it Coboundary expanders}, Journal of Topology and Analysis, {\bf 4}(2012), No.4, 499-514.
  \bibitem{GS}
    A Gundert and M Szedl\'ak, {\it Higher Dimensional Cheeger Inequalities}, Annual Symposium on Computational Geometry (NewYork, NY, USA), SOCG14, ACM 2014.
  \bibitem{Gro}
    M. Gromov, {\it Singularities, expanders and topology of maps. part 2: From combinatorics
to topology via algebraic isoperimetry}, Geometric and Functional Analysis {\bf 20}(2010), no. 2, 416-526.
  \bibitem{HLW}
    S. Hoory, N. Linial and A. Wigderson, {\it Expander graphs and there applications}, Bulletin of the American Mathematical Society {\bf 43}(4)(2008) 439-561.
  \bibitem{JLYY}
    X.Jiang, L-H. Lim, Y. Yao and Y. Ye, {\it Statistical ranking and combinatorial Hodge theory}, Mathematical Programming {\bf 127}(2011), 203-244.
  \bibitem{LM}
     N.Linial and R. Meshulam, {\it Homological connectivity of random 2-complexes}, Combinatorica {\bf 26}(2006), no.4, 475-487. 
  \bibitem{L}
    A. Lubotzky, {\it Ramanujan complexes and high dimensional expanders}, Japaniese Journal of Mathematics {\bf 9}(2014), 137-169.
  \bibitem{PRT}
    O. Parzanchevski, R. Rosenthal and  R.J. Tessler, {\it Isoperimetric inequalities in simplicial complexes}, Combinatorica {\bf 36}(2) (2016) 195-227.
   \bibitem{SKM}
      J. Steenbergen, C. Klivans, and S. Mukherjee, {\it A Cheeger-type inequality on simplicial
complexes}, Advances in Applied Mathematics {\bf 56} (2014) 56-77.
  \bibitem{Ta}
	R.M. Tanner, {\it Explicit concentrators from generalized n-gons}, SIAM Journal on Algebraic
and Discrete Methods {\bf 5} (1984), 287.


\end{thebibliography}
\end{document}